\documentclass[11pt]{amsart}
\usepackage{amsmath,amssymb,latexsym,cite}
\usepackage{verbatim}
\usepackage[left=3.1cm,right=3.1cm,top=3.2cm,bottom=3.2cm]{geometry}

\usepackage{color,graphicx}
\usepackage[colorlinks=true,urlcolor=blue, citecolor=red,linkcolor=blue,linktocpage,pdfpagelabels, bookmarksnumbered,bookmarksopen]{hyperref}
\usepackage[hyperpageref]{backref}
\usepackage[english]{babel}
\makeatletter
\providecommand\@dotsep{5}
\def\listtodoname{List of Todos}
\def\listoftodos{\@starttoc{tdo}\listtodoname}
\makeatother

\renewcommand{\phi}{{\varphi}}

\renewcommand{\epsilon}{{\varepsilon}}
\renewcommand{\theta}{{\vartheta}}

\numberwithin{equation}{section}

\newtheorem{theorem}{Theorem}[section]
\newtheorem{lemma}[theorem]{Lemma}

\newtheorem{remark}[theorem]{Remark}

\begin{document}

\title[Exponential nonlinearity  combined with convection term]{Positive solutions of quasilinear elliptic equations with exponential nonlinearity  combined with convection term}


\author{Anderson L.A. de Araujo}
\address{Universidade Federal de Vi\c{c}osa, Departamento de Matem\'atica, Avenida Peter Henry Rolfs, s/n, CEP 36570-900, Vi\c{c}osa, MG, Brazil}
\curraddr{}
\email{anderson.araujo@ufv.br}
\thanks{}

\author{Luiz F.O Faria}
\address{Departamento de Matem\'atica, Universidade Federal de Juiz de Fora, CEP 30161-970, Juiz de Fora - MG, Brazil}
\curraddr{}
\email{luiz.faria@ufjf.edu.br}
\thanks{Luiz F.O Faria was partially supported by FAPEMIG CEX APQ 02374/17.}



\date{}

\dedicatory{}

\begin{abstract}
We establish the existence of positive solutions for a nonlinear elliptic Dirichlet problem in dimension $N$ involving the $N$-Laplacian. The nonlinearity considered depends on the gradient of the unknown function and an exponential term. In such case, variational methods cannot be applied. Our approach is based on approximation scheme, where we consider a new class of normed spaces of finite dimension. As a particular case, we extended the result achieved by De Araujo and Montenegro [2016]  for any  $N>2$.
\end{abstract}

\maketitle


\section{Introduction}

Let $\Omega \subset \mathbb{R}^N$ be a smooth bounded domain and $p>1$. Consider the following problem

\begin{equation}\label{m1}\left\{
    \begin{array}{lll}
        -\Delta_pu=g(x,u,\nabla u) &\mbox{in}&  \Omega, \\
        u=0 & \mbox{on} & \partial\Omega.
    \end{array} \right.
\end{equation}
Here, the  operator $-\Delta_p:W^{1,p}_0(\Omega)\to
W^{-1,p'}(\Omega)$,  where $\frac{1}{p}+\frac{1}{p'}=1$, is defined by
\begin{eqnarray*}
\langle-\Delta_pu,v\rangle=\int_\Omega|\nabla u|^{p-2}\nabla u\nabla
v dx\, \mbox{ for all $u,v\in W^{1,p}_0(\Omega)$},
\end{eqnarray*}
and the forcing term $g$ has the form of a convection term, that is, it depends also on the gradient of the unknown function. Due to the presence of the gradient $ \nabla u$ in the term $g(x, u,  \nabla u)$, problem \eqref{m1} does not have, in general, variational structure. This kind of problems
are usually studied by means of topological degree, the method of sub-supersolutions, fixed point
theory, approximation techniques and iterative scheme. For instance, we would like to cite \cite{amann-crandall, Brezis-g, Ruiz,fmm,Fig-gir-mat,Faraci-Mot-Pugl,FMMT}. In particular, in \cite{fmm}, via an approximation on finite dimensional subspaces,  the authors proved the existence of a positive solution for the following problem
\begin{equation*}\left\{
    \begin{array}{lll}
        -\Delta_pu-\mu\Delta_qu=g(x,u,\nabla u) &\mbox{in}&  \Omega,\nonumber\\
        u>0&\mbox{in}&  \Omega,\nonumber\\
        u=0 & \mbox{on} & \partial\Omega,\nonumber
    \end{array} \right.\leqno(P)
\end{equation*}
where $\mu\geq 0$, $g:\Omega\times \mathbb{R}\times \mathbb{R}^{N}\to
\mathbb{R}$ is a continuous function satisfying the growth
condition:
    \begin{itemize}
                    \item[$(G)$]  $b_0|t|^{r_0}\leq f(x,t,\xi)\leq b_{1} (1+|t|^{r_1}+|\xi|^{r_2})$ \;
                    for all $(x,t,\xi)\in\Omega\times \mathbb{R}\times\mathbb{R}^N$,
                    with constants $b_0,b_1>0$, $r_{1},r_2\in [0,p-1)$, $r_0\in [0,p-1)$
                    if $\mu=0$,
                    and $r_0\in [0,q-1]$ if $\mu>0$.
    \end{itemize}


On the other hand, elliptic problems of the type
\begin{equation}\label{Pgeral}
	\left\{
	\begin{array}{lcc}
	-\Delta_N v = g(x,v) &\textup{in}& \Omega, \\
		v=0 &\textup{on}& \partial\Omega,
		\end{array}
		\right.
\end{equation}
where $\Omega\subset \mathbb{R}^N$ and  $g(x,v)$ is continuous and behaves like $\exp(\alpha |v|^{N/(N-1)})$ as $|v| \rightarrow +\infty$ have been studied by many authors, we would like to cite \cite{ap,FMR,FOR,F,LL,SS,AM}. One of the main ingredients is the Trudinger-Moser inequality introduced in \cite{m,t}. Namely, given $u \in W^{1,N}_0(\Omega)$, then
\begin{equation}\label{MT1}
e^{\sigma|u|^{\frac{N}{N-1}}} \in L^1(\Omega)\,\,\,\mbox{for every}\,\,\,\sigma >0,	
\end{equation}
and there exists a positive constant $L(N)$ which depends on $N$ only, such that
 \begin{equation}\label{MT2}
\sup_{\|u\|_{W^{1,N}_0(\Omega)}\leq 1}\int_{\Omega}e^{\sigma|u|^{\frac{N}{N-1}}}dx \leq L(N)|\Omega|\,\,\,\mbox{for every}\,\,\,\sigma\leq \alpha_N,
\end{equation}
where $|\Omega|=\int_{\Omega}dx$, $\alpha_N=Nw_{N-1}^{\frac{1}{N-1}}$ and $w_{N-1}$ is the $(N-1)$-dimensional measure of the $(N-1)$-sphere.

 In particular, in \cite{AM} the authors proved existence of solutions for the following problem

\begin{equation*}\left\{
    \begin{array}{lll}
        -\Delta u=\lambda u^q+f(u) &\mbox{in}&  \Omega,\nonumber\\
        u>0&\mbox{in}&  \Omega,\nonumber\\
        u=0 & \mbox{on} & \partial\Omega,\nonumber
    \end{array} \right.
\end{equation*}
where $\Omega\subset\mathbb{R}^2$, $\lambda>0$ is a small enough parameter, $0<q<1$ and
 $f: [0,\infty)\to
\mathbb{R}$ is a continuous function satisfying the growth
condition:
    \begin{itemize}
                    \item[$(H)$] $0\leq tf(t)\leq C|t|^{r}\exp(\alpha t^{2})$ \
                    where
                     $\alpha>0$ and $r>2$.
  \end{itemize}
In \cite{AM2}, still considering $N=2$, the authors proved existence of solutions for an elliptic system with arguments based in \cite{AM} and with nonlinearities satisfying the growth condition ($H$).

In this work we are concerned with the existence of positive solutions for the
problem:
\begin{equation*}\left\{
    \begin{array}{lll}
        -\Delta_Nu=\lambda (a_1 u^{r_1}+a_2|\nabla u|^{r_2})+f(u) &\mbox{in}&  \Omega\nonumber\\
         u=0 & \mbox{on} & \partial\Omega,\nonumber
    \end{array} \right.\leqno(P)
\end{equation*}
where $\Omega$ in $\mathbb{R}^N$ is a bounded domain with a
$C^{1,\alpha}$-boundary $\partial\Omega$, for some $0<\alpha\leq 1$, $\lambda>0$ is a parameter, $0<r_i<N-1$, for $i=1,2$,  $a_1>0$, $a_2\geq 0$, and
 $f: [0,\infty)\to
\mathbb{R}$ is a nonegative continuous function. The main assumption on the function $f$ is the following, which will be referred
throughout the paper as $(F)$:\\ 
    \begin{itemize}
                    \item[$(F)$] $0\leq tf(t)\leq a_3t^{r_3+1}\exp(\alpha t^{\frac{N}{N-1}})$ \
                    where
                     $a_3,\alpha>0$, and $r_3>N-1$.
  \end{itemize}

 Most of the papers, to prove existence results for the problem, assume  Ambrosetti--Rabinowitz conditions (or some additional conditions) to obtain Palais-Smale or Cerami compactness condition. Notice that in this paper we don't need to impose such extra hypotheses.

An interesting problem related to (P), by considering a more general operator, was treated by \cite{YP}. The authors studied a $(N,q)$-Laplacian problem with a critical Trundinger-Moser nonlinearity as following
 \begin{equation*}\left\{
    \begin{array}{lll}
        -\Delta_Nu-\Delta_q u=\mu  |u|^{q-2}u+\lambda|u|^{N-2}ue^{|u|^{N/(N-1)}} &\mbox{in}&  \Omega\nonumber\\
         u=0 & \mbox{on} & \partial\Omega,\nonumber
    \end{array} \right.
\end{equation*}
where $N>q>1$,  $\mu \in \mathbb{R}$ and $\lambda>0$. By using a critical point theorem, based on a cohomological index, they proved  the existence of solution for $\mu $ interacting with the  first eigenvalue of the $(-\Delta_q u,W_0^{1,q}(\Omega))$ operator and for $\lambda$ sufficiently large.

Here we extend  the results of \cite{AM} for the general dimension case $N>1$ ($a_2=0$). In order to prove the existence of positive solutions for $(P)$, we borrow some ideas from \cite{AM} and \cite{fmm}. Due to the presence of the supercritical term $\exp(\alpha |v|^{N/(N-1)})$, along  with the convection term, we had to overcome some problems. For example, in $W_0^{1,N}(\Omega)$ we need to assume a Schauder basis instead of the Hilbert basis (like in \cite{AM}), which becomes some additional difficulty. By comparing with \cite{fmm}, due to the presence of the term $\exp(\alpha |v|^{N/(N-1)})$, a suitable modification on the approximating approach had to be done. Although in  \cite{fmm} the authors used the Schauder basis, we could not obtain the necessary estimates for this approach by considering the approximate spaces used there. To do this, we consider a new class of normed spaces of finite dimension. 

Our main result reads as follows:

\begin{theorem}\label{TP}
 Suppose that $f:[0,\infty)\to
\mathbb{R}$ is a continuous function satisfying the  assumption $(F)$. Then there exists $\lambda^{*}>0$ such that for every $\lambda\in (0,\lambda^{*})$ problem $(P)$ admits a (positive) weak solution
$u\in W^{1,N}_0(\Omega)$.
\end{theorem}

\section{Preliminary results} \label{s.2}
The Sobolev space $W^{1,N}_0(\Omega)$ is endowed with the norm

$$\|u\|_{W^{1,N}_0(\Omega)}=\left(\int _{\Omega}|\nabla u|^N dx \right)^{1/N}.$$
To prove Theorem \ref{TP} we approximate $f$ by Lipschitz functions $f_k:\mathbb{R} \to \mathbb{R}$ defined by
\begin{equation}\label{eq111}
f_k(s)=\displaystyle\left\{
	\begin{array}{lcc}
-k[G(-k-\frac{1}{k}) - G(-k)], &\textup{if}& s\leq -k,\\	
-k[G(s-\frac{1}{k}) - G(s)], &\textup{if}& -k\leq s \leq -\frac{1}{k},\\
k^2s[G(-\frac{2}{k}) - G(-\frac{1}{k})], &\textup{if}& -\frac{1}{k}\leq s\leq 0,\\
k^2s[G(\frac{2}{k}) - G(\frac{1}{k})], &\textup{if}& 0\leq s\leq \frac{1}{k},\\	
k[G(s+\frac{1}{k}) - G(s)], &\textup{if}& \frac{1}{k}\leq s \leq k,\\
k[G(k+\frac{1}{k}) - G(k)], &\textup{if}& s\geq k,\\
	\end{array}
		\right.
\end{equation}
where $G(s)=\int_0^sf(\xi)d\xi$.

The following (approximation) result was proved in \cite{s} and uses the explicit expression of the sequence \eqref{eq111}.

\begin{lemma}\label{lemma1}
Let $f:\mathbb{R} \to \mathbb{R}$ be a continuous function such that $sf(s)\geq 0$ for every $s \in \mathbb{R}$. Then there exists a sequence $f_k:\mathbb{R} \to \mathbb{R}$ of continuous functions satisfying

(i) $sf_k(s)\geq 0$ for every $s \in \mathbb{R}$;

(ii) $\forall \, k \in \mathbb{N}$ $\exists c_k>0$ such that $|f_k(\xi) - f_k(\eta)|\leq c_k|\xi - \eta|$ for every $\xi, \eta \in \mathbb{R}$;
	
(iii) $f_k$ converges uniformly to $f$ in bounded subsets of $\mathbb{R}$.
\end{lemma}

The sequence $f_k$ of the previous lemma has some additional properties.
\begin{lemma}\label{lemma2}
Let $f: \mathbb{R} \to \mathbb{R}$ be a continuous function satisfying $(F)$ for every $s \in \mathbb{R}$.
Then the sequence $f_k$ of Lemma \ref{lemma1} satisfies

(i) $\forall \, k \in \mathbb{N}$, $0\leq sf_k(s) \leq C_1|s|^{ r_3}\exp(2^{\frac{N}{N-1}}\alpha\,|s|^{\frac{N}{N-1}})$ for every $|s|\geq \frac{1}{k}$;
	
(ii) $\forall \, k \in \mathbb{N}$, $0\leq sf_k(s) \leq C_2|s|^2\exp(2^{\frac{N}{N-1}}\alpha\,|s|^{\frac{N}{N-1}})$ for every $|s|\leq \frac{1}{k}$,

\noindent where $C_1$ and $C_2$ are positive constants independent of $k$.
\end{lemma}
\proof  Everywhere in this proof the constant $a_3$ is  the one of (\ref{eq111}).

\textit{\underline{First step}}. Suppose that $-k\leq s \leq -\frac{1}{k}$.

By the mean value theorem, there exists $\eta\in (s-\frac{1}{k},s)$ such that
\[
f_k(s)=-k[G(s-\frac{1}{k}) - G(s)]=-kG'(\eta)(s-\frac{1}{k}-s)=f(\eta)
\]
and
\[
sf_k(s)=sf(\eta).
\]
Since $s-\frac{1}{k}<\eta<s<0$ and $f(\eta)<0$, we have $sf(\eta) \leq \eta f(\eta)$. Therefore,
\[
\begin{array}{rcl}
sf_k(s) \leq \eta f(\eta) &\leq & a_3|\eta|^{r_3}\exp(\alpha\,|\eta|^{\frac{N}{N-1}})\\
 &\leq & a_3|s-\frac{1}{k}|^{r_3}\exp(\alpha\,|s-\frac{1}{k}|^{\frac{N}{N-1}})\\
&\leq & a_3(|s| +\frac{1}{k})^{r_3}\exp(\alpha\,(|s| +\frac{1}{k})^{\frac{N}{N-1}})\\
&\leq& a_3(2|s|)^{r_3}\exp(\alpha\,(2|s|)^{\frac{N}{N-1}})\\
&=& a_32^{r_3}|s|^{r_3}\exp(2^{\frac{N}{N-1}}\alpha\,|s|^{\frac{N}{N-1}}).
\end{array}
\]

\textit{\underline{Second step}}. Assume $\frac{1}{k}\leq s \leq k$.

By the mean value theorem, there exists $\eta\in (s,s+\frac{1}{k})$ such that
\[f_k(s)=k[G(s+\frac{1}{k}) - G(s)]=kG'(\eta)(s+\frac{1}{k}-s)=f(\eta)\]
and
\[sf_k(s)=sf(\eta).\]
Since $0<s <\eta <s+\frac{1}{k}$ and $f(\eta)>0$, we have $sf(\eta) \leq \eta f(\eta)$. Therefore,
\[
\begin{array}{rcl}
sf_k(s) \leq \eta f(\eta) &\leq & a_3|\eta|^{r_3}\exp(\alpha\,|\eta|^{\frac{N}{N-1}})\\
 &\leq & a_3|s+\frac{1}{k}|^{r_3}\exp(\alpha\,|s+\frac{1}{k}|^{\frac{N}{N-1}})\\
&\leq& a_3(2|s|)^{r_3}\exp(\alpha\,(2|s|)^{\frac{N}{N-1}})\\
&=& a_32^{r_3}|s|^{r_3}\exp(2^{\frac{N}{N-1}}\alpha\,|s|^{\frac{N}{N-1}}).
\end{array}
\]

\textit{\underline{Third step}}. Suppose that $|s|\geq k$, then
\begin{equation}\label{eq3}
	f_k(s)=\displaystyle\left\{
	\begin{array}{lcc}
	-k[G(-k-\frac{1}{k}) - G(-k)], &\textup{if}& s\leq -k\\
	k[G(k+\frac{1}{k}) - G(k)], &\textup{if}& s\geq k.\\
		\end{array}
		\right.
\end{equation}

If $s\leq -k$, by the mean value theorem, there exists $\eta\in (-k-\frac{1}{k},-k)$ such that
\[f_k(s)=k[G(-k-\frac{1}{k}) - G(-k)]=-kG'(\eta)(-k-\frac{1}{k}-(-k))=f(\eta)\]
and
\[sf_k(s)=sf(\eta).\]
Since $-k-\frac{1}{k} <\eta <-k<0$ and $k<|\eta| < k + \frac{1}{k}$, we conclude that
\begin{equation}\label{eq3.1}
\begin{array}{rcl}
sf_k(s)=\frac{s}{\eta}\eta f(\eta) &\leq& \frac{|s|}{|\eta|}a_3|\eta|^{r_3}\exp(\alpha\,|\eta|^{\frac{N}{N-1}}) =a_3|s||\eta|^{{r_3}-1}\exp(\alpha\,|\eta|^{\frac{N}{N-1}})\\
 &\leq& a_3|s|(k + \frac{1}{k})^{{r_3}-1}\exp(\alpha\,(k + \frac{1}{k})^{\frac{N}{N-1}})\\
 &\leq& a_3|s|(|s| + \frac{1}{k})^{{r_3}-1}\exp(\alpha\,(|s| + \frac{1}{k})^{\frac{N}{N-1}})\\
&\leq& a_3|s|(2|s|)^{{r_3}-1}\exp(\alpha\,(2|s|)^{\frac{N}{N-1}})\\
&\leq& a_32^{{r_3}-1}|s|^{r_3}\exp(2^{\frac{N}{N-1}}\alpha\,|s|^{\frac{N}{N-1}}).
\end{array}
\end{equation}

If $s\geq k$, by the mean value theorem, there exists $\eta\in (k,k + \frac{1}{k})$ such that
\[f_k(s)=k[G(k+\frac{1}{k}) - G(k)]=kG'(\eta)(k+\frac{1}{k}-k)=f(\eta).\]
By computations similar to conclude (\ref{eq3.1}) one has
\[sf_k(s)=sf(\eta)=\frac{s}{\eta}\eta f(\eta) \leq \frac{|s|}{|\eta|}a_3|\eta|^{r_3} \exp(\alpha\,|\eta|^{\frac{N}{N-1}}) \leq a_32^{{r_3}-1}|s|^{r_3}\exp(2^{\frac{N}{N-1}}\alpha\,|s|^{\frac{N}{N-1}}).\]

\textit{\underline{Fourth step}}. Assume $-\frac{1}{k}\leq s\leq \frac{1}{k}$, then
\begin{equation}
	f_k(s)=\displaystyle\left\{
	\begin{array}{lcc}
	k^2s[G(-\frac{2}{k}) - G(-\frac{1}{k})], &\textup{if}& -\frac{1}{k}\leq s\leq 0\\
	k^2s[G(\frac{2}{k}) - G(\frac{1}{k})], &\textup{if}& 0\leq s\geq \frac{1}{k}.\\
		\end{array}
		\right.
\end{equation}

If $-\frac{1}{k}\leq s\leq 0$, by the mean value theorem, there exists $\eta\in (-\frac{2}{k},-\frac{1}{k})$ such that
\[f_k(s)=k^2s[G(-\frac{2}{k}) - G(-\frac{1}{k})]=k^2sG'(\eta)(-\frac{2}{k}-(-\frac{1}{k}))=-ksf(\eta).\]
Therefore
\[
sf_k(s)=-ks^2f(\eta)=-k\frac{s^2}{\eta}\eta f(\eta) \leq k\frac{s^2}{|\eta|}\eta f(\eta) \leq
\]
\[
\leq a_3k|s|^2|\eta|^{{r_3}-1}\exp(\alpha\,|\eta|^{\frac{N}{N-1}}) \leq a_3k|s|^2(\frac{2}{k})^{{r_3}-1}\exp(\alpha\,|\eta|^{\frac{N}{N-1}}) \leq a_32^{{r_3}-1}|s|^2\exp\left(\alpha\,\left(\frac{2}{k}\right)^{\frac{N}{N-1}}\right)
\]
\begin{equation}\label{eq3.2}
\leq a_32^{{r_3}-1}|s|^2\exp(2^{\frac{N}{N-1}}\alpha) \leq a_32^{{r_3}-1}\exp(2^{\frac{N}{N-1}}\alpha)|s|^2\exp(2^{\frac{N}{N-1}}\alpha\,|s|^{\frac{N}{N-1}}).
\end{equation}

If $0\leq s \leq \frac{1}{k}$, by the mean value theorem, there exists $\eta\in (\frac{1}{k},\frac{2}{k})$ such that
\[f_k(s)=k^2s[G(\frac{2}{k}) - G(\frac{1}{k})]=k^2sG'(\eta)(\frac{2}{k}-\frac{1}{k})=ksf(\eta).\]
By similar computations to conclude (\ref{eq3.2}) one obtains
\[sf_k(s)=ks^2f(\eta)=k\frac{s^2}{|\eta|}\eta f(\eta)  \leq a_32^{{r_3}-1}\exp(2^{\frac{N}{N-1}}\alpha)|s|^2\exp(2^{\frac{N}{N-1}}\alpha\,|s|^{\frac{N}{N-1}}). \]
The proof of the lemma follows by taking $C_1=a_32^{r_3}$ ad $C_2=a_32^{{r_3}-1}C2^{{r_3}-1}\exp(2^{\frac{N}{N-1}}\alpha)$ where $a_3$ is given in $(F)$. \qed

Before concluding this section, we will enunciate a comparison principle due to Faria, Miyagaki and Motreanu \cite[Theorem 2.2]{fmm}.

Consider the Dirichlet problem
\begin{equation}\label{f1}
\left\{ \begin{array}{lll}
-\Delta_p u-\mu\Delta_q u=g(u)&\mbox{ in }& \Omega\\
u=0 & \mbox{ on } & \partial \Omega, \end{array} \right.
\end{equation}
where $1<q\leq p<+\infty$, $\mu\geq 0$ and
$g:\mathbb{R}\rightarrow\mathbb{R}$ is a continuous function.

We recall that $u_1\in W^{1,p}(\Omega)$ is a subsolution of problem
(\ref{f1}) if $u_1\geq 0$ a.e. on $\partial\Omega$ and
$$
\int_\Omega(|\nabla u_1|^{p-2}\nabla u_1\nabla\varphi+\mu|\nabla
u_1|^{q-2}\nabla u_1\nabla\varphi) dx\leq \int_\Omega g(u_1)\varphi
dx
$$
for all $\varphi\in W^{1,p}_0(\Omega)$ with $\varphi\geq 0$ a.e. in
$\Omega$, while $u_2\in W^{1,p}(\Omega)$ is a supersolution of
(\ref{f1}) if the reversed inequalities are satisfied with $u_2$ in
place of $u_1$ for all $\varphi\in W^{1,p}_0(\Omega)$ with
$\varphi\geq 0$ a.e. in $\Omega$.

\begin{theorem}\label{teorsubsup}
Let $g:\mathbb{R}\rightarrow\mathbb{R}$ be a continuous function
such that $t^{1-q}g(t)$ is nonincreasing for $t>0$ if $\mu>0$, and
$t^{1-p}g(t)$ is nonincreasing for $t>0$ if $\mu=0$. Assume that
$u_1\in W^{1,p}_0(\Omega)$ and $u_2\in W^{1,p}_0(\Omega)$ are a
positive subsolution and a positive supersolution of problem
(\ref{f1}), respectively. If $u_i\in L^{\infty}(\Omega)\cap
C^{1,\alpha}(\Omega)$, $\Delta_p u_i\in L^{\infty}(\Omega)$,
$u_i/u_j\in L^{\infty}(\Omega)$ for $i,j=1,2$, then $u_2\geq u_1$ in
$\Omega$.
\end{theorem}

\section{Approximation problem} \label{s.3}

For each $n\in \mathbb{N}$, we define the auxiliary problem $(P_n)$ by
\begin{equation*}
	\left\{
    \begin{array}{lll}
        -\Delta_Nu=\lambda (a_1 (u_+)^{r_1}+a_2|\nabla u|^{r_2})+f_n(u)+\frac{1}{n} &\mbox{in}&  \Omega\nonumber\\
        u>0&\mbox{in}&  \Omega\nonumber\\
        u=0 & \mbox{on} & \partial\Omega,\nonumber
    \end{array}
		\right. \leqno(P_n)
\end{equation*}
where $f_n$ are given by Lemma \ref{lemma1} and Lemma \ref{lemma2}, and $u_+ = \max\{u,0\}$.

To prove Theorem \ref{TP} we first show the existence of a solution of  problem $(P_n)$ by using the Galerkin method. We would like to cite \cite{af1} as the seminal paper in this type of approach.

\subsection{Finite-Dimensional Spaces}
Let $\mathcal{B}=\{w_1,w_2,\dots,w_n,\dots\}$ be a Schauder basis of $W^{1,N}_0(\Omega)$ (see \cite{FJN,Brezis}). For each positive integer $m$, let 
\[
W_m=[w_1,w_2,\dots,w_m]
\]  be the $m$-dimensional subspace of $W^{1,N}_0(\Omega)$ (generated by $\{w_1,w_2,\dots,w_m\}$)  with norm  induced from  $W^{1,N}_0(\Omega)$.  
Let $ \xi=(\xi_{1},\ldots,\xi_{m})\in \mathbb{R}^m$, notice that 
$$|\xi|_m = \|\sum _{j=1}^{m} \xi_jw_j\|_{W^{1,N}_0(\Omega)},$$
defines a norm in $\mathbb{R}^m$. In fact, let $\xi^i=(\xi^i_{1},\ldots,\xi^i_{m})\in \mathbb{R}^m$, $i=1,2$, and let $\lambda \in \mathbb{R}$.
\begin{itemize}
\item [(i)] $|\xi^1+\xi^2|_m \leq |\xi^1|_m+|\xi^2|_m$:
\begin{eqnarray*}
|\xi^1+\xi^2|_m &=& \|\sum _{j=1}^{m} \xi^1_jw_j+\sum _{j=1}^{m} \xi^2_jw_j\|_{W^{1,N}_0(\Omega)}\\ &\leq& \|\sum _{j=1}^{m} \xi^1_jw_j\|_{W^{1,N}_0(\Omega)}+\|\sum _{j=1}^{m} \xi^2_jw_j\|_{W^{1,N}_0(\Omega)} \\ &=& |\xi^1|_m+|\xi^2|_m.
\end{eqnarray*}

\item [(ii)] $|\lambda \xi^1|_m = |\lambda||\xi^1|_m$:
\begin{eqnarray*}
|\lambda\xi^1|_m = \|\lambda\sum _{j=1}^{m} \xi^1_jw_j\|_{W^{1,N}_0(\Omega)}=|\lambda|\|\sum _{j=1}^{m} \xi^1_jw_j\|_{W^{1,N}_0(\Omega)}= |\lambda||\xi^1|_m.
\end{eqnarray*}

\item [(iii)] $| \xi^1|_m = 0 \Leftrightarrow \xi^1=0$:\\

$(\Rightarrow)$ $0=| \xi^1|_m=\|\sum _{j=1}^{m} \xi^1_jw_j\|_{W^{1,N}_0(\Omega)} $ implies $\sum _{j=1}^{m} \xi^1_jw_j=0$. By uniqueness of the representation (using a Schauder basis) of the null vector,  we conclude that $\xi^1=0$. \\

$(\Leftarrow)$ It is trivial.

\end{itemize}

By using the above notation, we can identify the normed spaces $(W_m, \|\cdot \|_{W^{1,N}_0(\Omega)})$ and $(\mathbb{R}^m,|\cdot|_m)$ by the isometric linear
transformation
\begin{equation}\label{transf}
v=\sum_{j=1}^{m}\xi_{j}w_{j}\in V_m\mapsto
\xi=(\xi_{1},\ldots,\xi_{m})\in\mathbb{R}^{m}.
\end{equation}

 The lemma below is a consequence of Brouwer’s Fixed Point Theorem and its proof can be found in Kesavan \cite{Kesavan}.

\begin{lemma}\label{prop1}
Let $F: \mathbb{R}^d \rightarrow \mathbb{R}^d$ be a continuous function such that $\left\langle F(\xi),\xi\right\rangle\geq 0$ for every $\xi \in \mathbb{R}^d$ with $|\xi|=r$ for some $r>0$. Then, there exists $z_0$ in the closed ball $\overline{B}_r(0)$ such that $F(z_0)=0$.
\end{lemma}

\subsection{Existence}
The following result is concerning the existence result for the auxiliary problem $(P_n)$.
\begin{lemma}\label{teo aux}
There exists $\lambda^*>0$ and $n^* \in \mathbb{N}$ such that $(P_n)$ admits a (positive) weak solution
$v\in W^{1,N}_0(\Omega)\cap C^{1,\alpha}(\overline{\Omega})$, for some $0<\alpha<1$, for every $\lambda \in (0,\lambda^*)$ and $n\geq n^*$.
\end{lemma}
\proof  
Let $\mathcal{B}=\{w_1,w_2,\dots,w_n,\dots\}$ be a Schauder basis of $W^{1,N}_0(\Omega)$.  For each positive integer $m$, let $W_m=[w_1,w_2,\dots,w_m]$.  By using the isometric linear transformation \eqref{transf}, define the function $F:\mathbb{R}^m \to \mathbb{R}^m$ such that $F(\xi)=(F_1(\xi),F_2(\xi),\dots, F_m(\xi))$, 
where
 \begin{equation*}
\begin{array}{lll}
F_j(\xi)&=&\displaystyle\int_{\Omega}|\nabla u|^{N-2}\nabla u\nabla w_jdx - \lambda\left(a_1\int_{\Omega}(u_+)^{r_1}w_jdx +a_2\int_{\Omega}|\nabla u|^{r_2}w_j\right)dx \\&& - \displaystyle \int_{\Omega}f_n(u_+)w_j - \frac{1}{n}\int_{\Omega}w_jdx, \,\,\,\, j=1,\ldots,m.
\end{array}
\end{equation*}
Therefore,
\begin{equation}\label{eq4}
\left\langle F(\xi),\xi \right\rangle=\int_{\Omega}|\nabla u|^{N}dx - \lambda\left(a_1\int_{\Omega}(u_+)^{r_1}udx +a_2\int_{\Omega}|\nabla u|^{r_2}udx\right)- \int_{\Omega}f_n(u_+)udx - \frac{1}{n}\int_{\Omega}udx.	
\end{equation}

Given $u \in W_m$, we define
$$
\Omega^+_n=\{x \in \Omega : |u(x)|\geq \frac{1}{n}\}
$$
and
$$
\Omega^-_n=\{x \in \Omega : |u(x)|< \frac{1}{n}\}.
$$
Thus, we rewrite (\ref{eq4}) as
\[\left\langle F(\xi),\xi\right\rangle = \left\langle F(\xi),\xi\right\rangle_P + \left\langle F(\xi),\xi\right\rangle_{N},\]
where
\[\left\langle F(\xi),\xi\right\rangle_P=\int_{\Omega^+_n}|\nabla u|^{N}dx - \lambda\left(a_1\int_{\Omega^+_n}(u_+)^{r_1}udx +a_2\int_{\Omega^+_n}|\nabla u|^{r_2}udx\right)- \int_{\Omega^+_n}f_n(u_+)u_+dx - \frac{1}{n}\int_{\Omega^+_n}udx\]
and
\[\left\langle F(\xi),\xi\right\rangle_{N}=\int_{\Omega^-_n}|\nabla u|^{N}dx - \lambda\left(a_1\int_{\Omega^-_n}(u_+)^{r_1}udx +a_2\int_{\Omega^-_n}|\nabla u|^{r_2}udx\right)- \int_{\Omega^-_n}f_n(u_+)u_+dx - \frac{1}{n}\int_{\Omega^-_n}udx.\]


\textit{\underline{Step 1}}. Since $0<r_i<N-1$, for $i=1,2$, then
\begin{equation}\label{eq4.1}
	\int_{\Omega^+_n}(u_+)^{r_1+1}dx\leq \int_{\Omega}(u_+)^{r_1+1}dx \leq \int_{\Omega}|u|^{r_1+1}dx = \|u\|^{r_1+1}_{L^{r_1+1}(\Omega)}\leq C_1\|u\|_{W^{1,N}_0(\Omega)}^{r_1+1}.
\end{equation}
By virtue of Lemma \ref{lemma2} (i) we get
\begin{equation}\label{eq4.2}
\begin{array}{rcl}
\displaystyle\int_{\Omega^+_n}f_n(u_+)u_+dx &\leq&\displaystyle C_1\int_{\Omega^+_n}|u|^{r_3+1}\exp(2^{\frac{N}{N-1}}\alpha|u|^{\frac{N}{N-1}})dx\\
&\leq&\displaystyle a_3 \left(\int_{\Omega}|u|^{N'(r_3+1)}dx\right)^{\frac{1}{N'}}\left(\int_{\Omega}\exp(N2^{\frac{N}{N-1}}\alpha|u|^{\frac{N}{N-1}})dx\right)^{\frac{1}{N}}\\
&=&\displaystyle a_3\|u\|_{L^{N'(r_3+1)}(\Omega)}^{r_3+1}\left(\int_{\Omega}\exp(N2^{\frac{N}{N-1}}\alpha|u|^{\frac{N}{N-1}})dx\right)^{\frac{1}{N}},
\end{array}
\end{equation}
where $\frac{1}{N}+\frac{1}{N'}=1$.

It follows from (\ref{eq4.1}) and (\ref{eq4.2}) that
\begin{equation}\label{eq12}
\begin{array}{rcl}
\displaystyle\left\langle F(\xi),\xi\right\rangle_P  &\geq& \displaystyle \int_{\Omega^+_n}|\nabla u|^{N}dx - \lambda\left(a_1C_1\|u\|_{W^{1,N}_0(\Omega)}^{r_1+1}+a_2\int_{\Omega^+_n}|\nabla u|^{r_2}udx\right) \\ &&- C_3\|u\|_{W^{1,N}_0(\Omega)}^{r_3+1}\left(\int_{\Omega}\exp(N2^{\frac{N}{N-1}}\alpha|u|^{\frac{N}{N-1}})dx\right)^{\frac{1}{N}} 
- \displaystyle \frac{C_4}{n}\|u\|_{W^{1,N}_0(\Omega)},
\end{array}
\end{equation}
where $C_0$, $C_1$ and $C_3$ are constants not depending $n$ and $m$.

\textit{\underline{Step 2}}. Since $0<r_i<N-1$, for $i=1,2$, then
\begin{equation}\label{eq10}
\int_{\Omega^-_n}(u_+)^{r_1+1} \leq \int_{\Omega^-_n}|u|^{r_1+1} \leq |\Omega|\frac{1}{n^{r_1+1}}.
\end{equation}
By virtue of Lemma \ref{lemma2} (ii) we get
\begin{equation}\label{eq11}
\int_{\Omega^-_n}f_n(u_+)u_+ \leq C_2 \int_{\Omega^-_n}|u|^2\exp(2^{\frac{N}{N-1}}\alpha|u|^{\frac{N}{N-1}})dx\leq C_2 \exp(2^{\frac{N}{N-1}}\alpha)|\Omega|\frac{1}{n^2}.
\end{equation}
It follows from (\ref{eq10}) and (\ref{eq11}) that
\begin{equation}\label{eq13}
	\left\langle F(\xi),\xi\right\rangle_N \geq \int_{\Omega^-_n}|\nabla u|^N -\lambda\left(a_1|\Omega|\frac{1}{n^{r_1+1}}+ a_2\int_{\Omega^-_n}|\nabla u|^{r_2}udx\right) \\ - C_3\exp(2^{\frac{N}{N-1}}\alpha)|\Omega|\frac{1}{n^2} 
- \displaystyle|\Omega|\frac{1}{n^2}.
\end{equation}
Since
\[
\int_{\Omega^+_n}|\nabla u|^{r_2}udx + \int_{\Omega^-_n}|\nabla u|^{r_2}udx=\int_{\Omega}|\nabla u|^{r_2}udx
\]
and
\begin{equation}\label{eq4.1.2}
\int_{\Omega}|\nabla u|^{r_2}|u|dx\leq \left(\int_{\Omega}|\nabla u|^{N}dx \right)^{r_2/N} \left(\int_{\Omega}|u|^{N/(N-r_2)}dx \right)^{(N-r_2)/N}\leq C \|u\|_{W^{1,N}_0(\Omega)}^{r_2+1}.
\end{equation}

Thus (\ref{eq12}), (\ref{eq13}) and (\ref{eq4.1.2}) imply
\begin{equation}\label{eq1}
\begin{array}{rcl}
\displaystyle\left\langle F(\xi),\xi\right\rangle  &\geq& \displaystyle \|u\|_{W^{1,N}_0(\Omega)}^{N} - \lambda\left(a_1C_1\|u\|_{W^{1,N}_0(\Omega)}^{r_1+1}+a_2C_2 \|u\|_{W^{1,N}_0(\Omega)}^{r_2+1}\right) \\ 
&-& \displaystyle C_3\|u\|_{W^{1,N}_0(\Omega)}^{r_3+1}\left(\int_{\Omega}\exp(N2^{\frac{N}{N-1}}\alpha|u|^{\frac{N}{N-1}})dx\right)^{\frac{1}{N}} 
- \displaystyle \frac{C_4}{n}\|u\|_{W^{1,N}_0(\Omega)}\\
&-& \displaystyle\lambda a_1|\Omega|\frac{1}{n^{r_1+1}} - C_5\exp(2^{\frac{N}{N-1}}\alpha)|\Omega|\frac{1}{n^2} 
- \displaystyle|\Omega|\frac{1}{n^2}.
\end{array}
\end{equation}

Assume now that $\|u\|_{W^{1,N}_0(\Omega)}=r$ for some $r>0$ to be chosen later. We have
\begin{equation}\label{eq7.2}
	\int_{\Omega}\exp(N2^{\frac{N}{N-1}}\alpha|u|^{\frac{N}{N-1}})dx = \int_{\Omega}\exp\left(N2^{\frac{N}{N-1}}\alpha\,r^{\frac{N}{N-1}}\left(\frac{|u|}{\|u\|_{W^{1,N}_0(\Omega)}}\right)^{\frac{N}{N-1}}\right)dx
\end{equation}
and in order to apply the Trudinger-Moser inequality (\ref{MT2}) we must have $N2^{\frac{N}{N-1}}\alpha\,r^{\frac{N}{N-1}} \leq \alpha_N$. Consequently,
\[
r\leq\frac{1}{2}\left(\frac{\alpha_N}{N\alpha}\right)^{\frac{N-1}{N}}.
\]
Then
\[
\sup_{\|u\|_{W^{1,N}_0(\Omega)}\leq 1}\int_{\Omega}\exp\left(N2^{\frac{N}{N-1}}\alpha\,r^{\frac{N}{N-1}}|u|^{\frac{N}{N-1}}\right)dx \leq L(N)|\Omega|.
\]

Hence,
\[
\begin{array}{rcl}
\left\langle F(\xi),\xi\right\rangle &\geq& \displaystyle r^N - \lambda(a_1C_1r^{r_1+1}+a_2C_2r^{r_2+1}) - C_3r^{r_3+1}L^{1/N}(N) - \displaystyle \frac{C_4}{n}r\\
&-& \displaystyle\lambda a_1|\Omega|\frac{1}{n^{r_1+1}} - C_5\exp(2^{\frac{N}{N-1}}\alpha)|\Omega|\frac{1}{n^2} 
- \displaystyle|\Omega|\frac{1}{n^2}.
\end{array}
\]
We need to choose $r$ such that
\[r^N - C_3 L(N)^{\frac{1}{N}}r^{r_3+1}\geq \frac{r^N}{2},\]
in other words,
\[
r\leq \frac{1}{(2 C_3 L(N)^{\frac{1}{N}})^{\frac{1}{r_3+1-N}}}.
\]
Thus, let $r=\min\left\{\frac{1}{2(2 C_3 L(N)^{\frac{1}{N}})^{\frac{1}{r_3+1-N}}}, \frac{1}{2}\left(\frac{\alpha_N}{N\alpha}\right)^{\frac{N-1}{N}}\right\}$, hence
\[
\left\langle F(\xi),\xi\right\rangle \geq \frac{r^N}{2} - \lambda(a_1C_1r^{r_1+1}+a_2C_2r^{r_2+1}) - \displaystyle \frac{C_4}{n}r - \displaystyle\lambda a_1|\Omega|\frac{1}{n^{r_1+1}} - C_5\exp(2^{\frac{N}{N-1}}\alpha)|\Omega|\frac{1}{n^2} 
- \displaystyle|\Omega|\frac{1}{n^2}.
\]
Now, defining  $\rho=\frac{r^N}{2} - \lambda(a_1C_1r^{r_1+1}+a_2C_2r^{r_2+1})$, we choose $\lambda^*>0$ such that $\rho>0$ for $\lambda < \lambda^*$.
Since $0<r_i<N-1$, for $i=1,2$, we can choose
$$
\lambda^*=\frac{1}{2}\frac{r^N}{2a_1C_1r^{r_1+1}+2a_2C_2r^{r_2+1}}.
$$
Now we choose $n^* \in \mathbb{N}$ such that
\[
\displaystyle \frac{C_4}{n}r + \displaystyle\lambda a_1|\Omega|\frac{1}{n^{r_1+1}} + C_5\exp(2^{\frac{N}{N-1}}\alpha)|\Omega|\frac{1}{n^2} 
+ \displaystyle|\Omega|\frac{1}{n^2} <\frac{\rho}{2},
\]
for every $n\geq n^*$. Let $\xi \in \mathbb{R}^m$, such that $|\xi|_m:=\left\|\sum_{i=1}^m\xi_iw_i\right\|_{W^{1,N}_0(\Omega)}=r$, then for $\lambda < \lambda^*$ and  $n\geq n^*$ we obtain
\begin{equation}\label{eq8}
	\left\langle F(\xi),\xi\right\rangle \geq \frac{\rho}{2}>0.
\end{equation}

Then by Lemma \ref{prop1}, for every $m \in \mathbb{N}$ there exists $y \in \mathbb{R}^m$ (with $|y|_m\leq r$) such that $F(y)=0$. Therefore, there exists $u_m \in W_m$ verifying
\begin{equation}\label{bond}
\|u_m\|_{W^{1,N}_0(\Omega)}\leq r ,\,\,  \mbox{for every} \,\, m \in \mathbb{N},
\end{equation}
and such that
\begin{equation}\label{eq15}
\begin{array}{lll}
\displaystyle \int_{\Omega}|\nabla u_m|^{N-2}\nabla u_m\nabla w dx&=& \displaystyle\lambda\left(a_1\int_{\Omega}(u_m^+)^{r_1}wdx +a_2\int_{\Omega}|\nabla u_m|^{r_2}wdx\right) \\ && +\displaystyle\int_{\Omega}f_n((u_m)_+)wdx + \frac{1}{n}\int_{\Omega}wdx, \,\, \forall \, w \in W_m.
\end{array}
\end{equation}
Since $W_m \subset W^{1,N}_0(\Omega)$ $\forall \, m \in \mathbb{N}$ and $r$ does not depend on $m$, then $(u_m)$ is a bounded sequence in $W^{1,N}_0(\Omega)$. Then, for some subsequence, there exists $u_n \in W^{1,N}_0(\Omega)$ (to simplify the notation, until the end of this section we will omit the subscript $n$ of the variable $u$) such that 
 \begin{equation}\label{eq16}
u_m \rightharpoonup u \,\,\, \mbox{weakly in} \,\,\, W^{1,N}_0(\Omega)
\end{equation}
and
 \begin{equation}\label{eq16.1}
u_m \to u \,\,\, \mbox{in} \,\,\, L^N(\Omega)\,\,\, \mbox{and a.e. in}\,\,\,\Omega.
\end{equation}

Notice that
\begin{equation}\label{eq30}
\|u\|_{W^{1,N}_0(\Omega)} \leq \liminf_{m \to \infty}\|u_m\|_{W^{1,N}_0(\Omega)}\leq r, \,\, \forall\, n \in \mathbb{N},
\end{equation}
and $r$ does not depend on $n$.  We claim that
 \begin{equation}\label{eq16.2}
u_m \to u \,\,\, \mbox{in} \,\,\, W^{1,N}_0(\Omega).
\end{equation}
Using the fact that  $\mathcal{B}=\{w_1,w_2,\dots,w_n,\dots\}$ is a Schauder basis of $W^{1,N}_0(\Omega)$, for every $u \in W^{1,N}_0(\Omega)$  there exists a unique sequence $(\alpha_n)_{n\geq 1}$ in $\mathbb{R}$  such that $u=\sum _{j=1}^{\infty}\alpha_j w_j $, so 
\begin{equation}\label{166}
\psi_m:=\sum_{j=1}^m \alpha _j w_j \rightarrow u \,\, \mbox{ in } W_0^{1,N}(\Omega)\,\,\mbox{ as }m\rightarrow \infty.\end{equation}
Using as test function $w=(u_m-\psi_m)\in W_m$ in \eqref{eq15}, we get

\begin{equation}\label{eq117}\begin{array}{lll}
\displaystyle\int_{\Omega}|\nabla u_m|^{N-2}\nabla u_m\nabla (u_m-\psi_m)dx &=&\displaystyle \lambda\left(a_1\int_{\Omega}(u_m^+)^{r_1}(u_m-\psi_m)dx +a_2\int_{\Omega}|\nabla u_m|^{r_2}(u_m-\psi_m)dx\right)\\ && \displaystyle + \int_{\Omega}f_n((u_m)_+)(u_m-\psi_m)dx + \frac{1}{n}\int_{\Omega}(u_m-\psi_m)dx.
\end{array}
\end{equation}

By continuity of $f_n$, \eqref{eq16}, \eqref{eq16.1}, \eqref{166} and hypothesis $(F)$, we get

\begin{equation}\label{eq113}
\lim _{m\rightarrow \infty} \frac{1}{n}\int_{\Omega}(u_m-\psi_m)dx=0,
\end{equation}
\begin{equation}\label{eq114}
\lim _{m\rightarrow \infty} \displaystyle a_1\int_{\Omega}(u_m^+)^{r_1}(u_m-\psi_m)dx=0,
\end{equation}
\begin{equation}\label{eq115}
\lim _{m\rightarrow \infty} a_2\int_{\Omega}|\nabla u_m|^{r_2}(u_m-\psi_m)dx=0,
\end{equation}
and
\begin{equation}\label{eq116}
\lim _{m\rightarrow \infty} \int_{\Omega}f_n((u_m)_+)(u_m-\psi_m)dx=0.
\end{equation}
Notice that \eqref{eq113}, \eqref{eq114} and \eqref{eq115} are immediately. Let us verigy \eqref{eq116}. By continuity of $f_n$ and (\ref{eq16.1}) we obtain
 \begin{equation*}
f_n((u_{m})_+)^{N'} \to f_n(u_+)^{N'} \,\,\, \mbox{ a.e. in}\,\,\, \Omega
\end{equation*}
and by Lemma \ref{lemma1} and (\ref{bond}), we obtain
 \begin{equation*}
\begin{array}{rcl}
\displaystyle\int_{\Omega}f_n((u_{m})_+)^{N'}dx&\leq&\displaystyle c_n^{N'}\int_{\Omega}|u_m|^{N'}=\|u_m\|_{L^{N'}(\Omega)}^{N'} \leq c_n^{N'}C\|u_m\|^{N'}_{W^{1,N}_0(\Omega)}\leq c_n^{N'}Cr^{N'}.\\
\end{array}
\end{equation*}
Hence, \cite[Theorem 13.44]{EK} leads to 
 \begin{equation}\label{eq16.7}
f_n((u_{m})_+) \to f_n(u_+) \,\,\, \mbox{ weakly in } \,\,\, L^{N'}(\Omega).
\end{equation}
Applying \eqref{eq16.1}, \eqref{166} and \eqref{eq16.7}, we conclude that \eqref{eq116} holds.

By \eqref{bond} and  \eqref{eq16}, we obtain
\begin{equation}\label{eq70}
\lim_{m\rightarrow\infty} \displaystyle\int_{\Omega}|\nabla u_m|^{N-2}\nabla u_m\nabla (u-\psi_m)dx=0. 
\end{equation}
By $\eqref{eq113}-\eqref{eq116}$ and  \eqref{eq70}, we obtain
\begin{equation}\label{eq71}
\lim_{m\rightarrow\infty} \displaystyle\int_{\Omega}|\nabla u_m|^{N-2}\nabla u_m\nabla (u_m-u)dx=0. 
\end{equation}
Now it is sufficient to apply the $(S_+)-$ property of $-\Delta_p$
(see, e.g., \cite[Proposition 3.5.]{MonMonPapa}) for obtaining
(\ref{eq16.2}).

Let	$k \in \mathbb{N}$, then for every $m\geq k$ we obtain
\begin{equation*}\label{eq17}
\begin{array}{lll}
\displaystyle\int_{\Omega}|\nabla u_m|^{N-2}\nabla u_m\nabla w_k dx &=&\displaystyle \lambda\left(a_1\int_{\Omega}(u_m^+)^{r_1}w_kdx +a_2\int_{\Omega}|\nabla u_m|^{r_2}w_kdx\right)+ \int_{\Omega}f_n((u_m)_+)w_kdx\\
&&\displaystyle + \frac{1}{n}\int_{\Omega}w_kdx, \,\,\, \forall \, w_k \in W_k.
\end{array}
\end{equation*}

Lettin $m\rightarrow \infty$, on accout of \eqref{eq16.2} we arrive at

\begin{equation*}\label{eq20}
\begin{array}{lll}
\displaystyle\int_{\Omega}|\nabla u|^{N-2}\nabla u\nabla w_k dx &=&\displaystyle \lambda\left(a_1\int_{\Omega}(u^+)^{r_1}w_kdx +a_2\int_{\Omega}|\nabla u|^{r_2}w_kdx\right)+ \int_{\Omega}f_n(u_+)w_kdx\\
&&\displaystyle + \frac{1}{n}\int_{\Omega}w_kdx, \,\,\, \forall \, w_k \in W_k.
\end{array}
\end{equation*}
Since $[W_k]_{k \in \mathbb{N}}$ is dense in $W^{1,N}_0(\Omega)$ we conclude that
\begin{equation*}\label{eq21}
\begin{array}{lll}
\displaystyle\int_{\Omega}|\nabla u|^{N-2}\nabla u\nabla w dx &=& \displaystyle\lambda\left(a_1\int_{\Omega}(u^+)^{r_1}wdx +a_2\int_{\Omega}|\nabla u|^{r_2}wdx\right)+ \int_{\Omega}f_n(u_+)wdx\\
&&\displaystyle + \frac{1}{n}\int_{\Omega}wdx, \,\,\, \forall \, w \in W^{1,N}_0(\Omega).
\end{array}
\end{equation*}
Furthermore, $u\geq 0$ in $\Omega$. In fact, since $u_- \in W^{1,N}_0(\Omega)$ then from (\ref{eq21}) we obtain
\begin{equation*}\label{eq221}
\begin{array}{lll}
\displaystyle\int_{\Omega}|\nabla u|^{N-2}\nabla u\nabla u_- dx &=& \displaystyle\lambda\left(a_1\int_{\Omega}(u^+)^{r_1}u_-dx +a_2\int_{\Omega}|\nabla u|^{r_2}u_-dx\right)+ \int_{\Omega}f_n(u_+)u_-dx\\
&&\displaystyle + \frac{1}{n}\int_{\Omega}u_-dx.
\end{array}
\end{equation*}
Hence
\[
- \|u_-\|_{W^{1,N}_0(\Omega)}^N = \lambda\left(a_1\int_{\Omega}(u^+)^{r_1}u_-dx +a_2\int_{\Omega}|\nabla u|^{r_2}u_-dx\right)+ \int_{\Omega}f_n(u_+)u_-dx + \frac{1}{n}\int_{\Omega}u_-dx \geq 0,
\]
because $\int_{\Omega}f_n(u_+)u_-dx=0$. Then $u_- \equiv 0$ a.e. in $\Omega$.


The first inequality in hypothesis $(F)$ and the equation in
$(P_n)$ guarantee that $u\not=0$. Here the presence of
$\frac{1}{n}>0$ is needed. Next, we observe that hypothesis $(F)$
allows us to refer to \cite[Theorem 7.1]{LU} from which we infer
that $u\in L^\infty(\Omega)$. Furthermore, the regularity result up
to the boundary in \cite[Theorem 1]{lieberman} and \cite[p. 320]{L}
ensures that $u\in C^{1,\beta}(\overline{\Omega})$ with some
$\beta\in(0,1)$. We also note that we may apply the strong maximum
principle in \cite[Theorem 5.4.1]{PS}. We are thus in a
position to apply \cite[Theorem 5.4.1]{PS} concluding that $u>0$ in
$\Omega$ because we know that $u\geq 0$, $u\not=0$, thereby $u$ is a
solution of problem $(P_n)$. This completes the proof. \qed

\begin{remark}
To apply \cite[Theorem 7.1]{LU} and infer
that $u\in L^\infty(\Omega)$, notice that it is necessary to consider the approximating functions  $f_n$ (given by Lemma \ref{lemma1}) instead of $f$. In fact, since
\[
F_n(v,p)=\lambda (a_1 v^{r_1}+a_2|p|^{r_2})+f_n(v)+\frac{1}{n}
\]
satisfies the inequality $(7.2)$ in \cite{LU}, because
\[
sign(v).F_n(v,p)\leq \lambda (a_1 |v|^{r_1}+a_2|p|^{r_2})+c_n|v|+\frac{1}{n}.
\] 
While $F(v,p)=\lambda (a_1 v^{r_1}+a_2|p|^{r_2})+f(v)$ does not necessarily satisfy such a hypothesis, in fact
\[
sign(v).F(v,p)\leq \lambda (a_1 |v|^{r_1}+a_2|p|^{r_2})+a_3|u|^{r_3}\exp(\alpha |v|^{\frac{N}{N-1}}).
\]
\end{remark}

\section{Proof of the main result} \label{s.4}

In this section we will prove Theorem \ref{TP}.  Consider the following problem 
\begin{equation}\label{P5}
	\left\{
    \begin{array}{lll}
        -\Delta_Nv =\lambda a_1 v^{r_1} &\mbox{in}&  \Omega\\
        v>0&\mbox{in}&  \Omega\\
        v=0 & \mbox{on} & \partial\Omega
    \end{array} \right.
\end{equation}
where $\lambda, a_1$ and $r_1$ were given in Theorem \ref{TP}. This problem admits a solution $v_0\in C^1_0(\overline{\Omega})$, see for instace \cite[Lemma 4.1]{fmm}. The function $v_0$ allows us to bound from below the solutions of $(P_n)$.

For each $\lambda\in(0,\lambda^*)$ and $n \in \mathbb{N}$ sufficiently large, by using Lemma \ref{teo aux}, we get that equation $(P_n)$ has a weak solution $u_n \in W^{1,N}_0(\Omega)\cap C^{1,\alpha}(\overline{\Omega})$ for some
$\alpha\in(0,1)$.

In view of \eqref{eq30}, we can argue as for \eqref{eq16.2}, to find a subsequence $n\rightarrow\infty$ such that the
corresponding sequence $\{u_n\}$ is strongly convergent:
\begin{equation}\label{eq16.2.2}
u_n \to u \,\,\, \mbox{in} \,\,\, W^{1,N}_0(\Omega).
\end{equation}

In fact,  for some subsequence, there exists $u \in W^{1,N}_0(\Omega)$ such that
 \begin{equation}\label{eq916}
u_n \rightharpoonup u \,\,\, \mbox{weakly in} \,\,\, W^{1,N}_0(\Omega)
\end{equation}
and, by Sobolev embedding for $1 \leq s < +\infty$,
\[
u_n \rightarrow u \,\,\mbox{in}\,\,L^s(\Omega) \,\, \mbox{and} \,\, \mbox{a.e. in} \,\, \Omega.
\]

Notice that
\begin{equation}\label{P3}
	\begin{array}{c} 
	\displaystyle\int_{\Omega}|\nabla u_n|^{N-2}\nabla u_n\nabla w \geq \lambda\,a_1\int_{\Omega}(u_n)^{r_1}w , \forall w \in W^{1,N}_0(\Omega) \mbox{ with } w\geq0.\\
			\end{array}
		\end{equation}



Since
\[
u_n \rightarrow u\,\, \mbox{a.e. in}\,\, \Omega,
\]
we have
\begin{equation}\label{eq26.2}
f_n(u_n(x)) \rightarrow f(u(x))\,\, \mbox{a.e. in}\,\, \Omega,
\end{equation}
by the uniform convergence of Lemma \ref{lemma1} ($iii$).

By Lemma \ref{lemma2}
 \begin{equation*}
\begin{array}{rcl}
\displaystyle\int_{\Omega}f_n(u_{n})^{N'}dx&=&\displaystyle \int_{\Omega_n^+}f_n(u_{n})^{N'}dx + \int_{\Omega_n^-}f_n(u_{n})^{N'}dx,\\
\end{array}
\end{equation*}
\begin{equation*}
\begin{array}{rcl}
\displaystyle\int_{\Omega^+_n}f_n(u_n)^{N'}dx &\leq&\displaystyle C_1^{\frac{N}{N-1}}\int_{\Omega^+_n}|u_n|^{(r_3-1)\frac{N}{N-1}}\exp(\frac{N}{N-1}2^{\frac{N}{N-1}}\alpha|u_n|^{\frac{N}{N-1}})dx\\
&\leq&\displaystyle C_1^{\frac{N}{N-1}} \left(\int_{\Omega}|u_n|^{(r_3-1)\frac{N}{N-2}}dx\right)^{\frac{N-2}{N-1}}\left(\int_{\Omega}\exp(N2^{\frac{N}{N-1}}\alpha|u_n|^{\frac{N}{N-1}})dx\right)^{\frac{1}{N-1}}\\
&=&\displaystyle C_1^{\frac{N}{N-1}}\|u_n\|_{L^{(r_3-1)\frac{N}{N-2}}(\Omega)}^{(r_3-1)\frac{N}{N-1}}\left(\int_{\Omega}\exp(N2^{\frac{N}{N-1}}\alpha|u_n|^{\frac{N}{N-1}})dx\right)^{\frac{1}{N-1}}\\
&\leq&\displaystyle C\|u_n\|_{_{W^{1,N}_0(\Omega)}}^{(r_3-1)\frac{N}{N-1}}\left(\int_{\Omega}\exp(N2^{\frac{N}{N-1}}\alpha|u_n|^{\frac{N}{N-1}})dx\right)^{\frac{1}{N-1}}
\end{array}
\end{equation*}
and
\begin{equation*}
\int_{\Omega^-_n}f_n(u_n)^{N'}dx \leq C_2^{\frac{N}{N-1}} \int_{\Omega^-_n}|u_n|^{\frac{N}{N-1}}\exp(\frac{N}{N-1}2^{\frac{N}{N-1}}\alpha|u_n|^{\frac{N}{N-1}})dx
\end{equation*}
\[
\leq C_2^{\frac{N}{N-1}} \exp(\frac{N}{N-1}2^{\frac{N}{N-1}}\alpha)|\Omega|\frac{1}{n^{\frac{N}{N-1}}}.
\]

Since $\|u_n\|_{_{W^{1,N}_0(\Omega)}}\leq r$, by the estimates before, we obtain
\[
\displaystyle\int_{\Omega}f_n(u_{n})^{N'}dx\leq C,
\]
for each $n$. Since $f_n(u_n(x)) \rightarrow f(u(x))\,\, \mbox{a.e. in}\,\, \Omega$, \cite[Theorem 13.44]{EK} leads to 
 \begin{equation}\label{eq25.2}
f_n(u_{n}) \to f(u) \,\,\, \mbox{ weakly in } \,\,\, L^{N'}(\Omega).
\end{equation}

Recall from (\ref{eq21}) that,  forall $ w \in W^{1,N}_0(\Omega)$, 
\begin{equation}\label{eq225}
\int_{\Omega}|\nabla u_n|^{N-2}\nabla u_n\nabla w = \lambda\left(a_1\int_{\Omega}(u_n)^{r_1}w +a_2\int_{\Omega}|\nabla u_n|^{r_2}w\right)+ \int_{\Omega}f_n(u_n)w + \frac{1}{n}\int_{\Omega}w.
\end{equation}
Taking $w=u_n - u$ in (\ref{eq225}), we obtain
\begin{equation}\label{eq71}
\lim_{m\rightarrow\infty} \displaystyle\int_{\Omega}|\nabla u_n|^{N-2}\nabla u_n\nabla (u_n-u)dx=0. 
\end{equation}
Now it is sufficient to apply the $(S_+)-$ property of $-\Delta_p$  for obtaining
\eqref{eq16.2.2}.

Then from \eqref{eq16.2.2} and \eqref{eq25.2} and the fact that $u_n$ solves $(P_n)$,  by passing to limit when $n \to +\infty$ we get that 

\begin{equation}\label{eq29}
\int_{\Omega}|\nabla u|^{N-2}\nabla u\nabla w = \lambda\left(a_1\int_{\Omega}(u)^{r_1}w +a_2\int_{\Omega}|\nabla u|^{r_2}w\right)+ \int_{\Omega}f(u)w, \,\,\, \forall \, w \in W^{1,N}_0(\Omega).
\end{equation}

Now, we are going to check that
$u>0$ in $\Omega$. Notice that, by (\ref{P5}) and (\ref{P3}), for each $n$ sufficiently large $u_n$ is a supersolution and $v_0$ is a subsolution of Problem \eqref{P5}. In order to apply Theorem \ref{teorsubsup}, we need to check
that $\frac{u_{n}}{v_0},
\frac{v_0}{u_{n}}\in L^\infty(\Omega)$. This follows by using Hopf boundary point lemma (in the strong maximum principle for
both Dirichlet problems \eqref{P5} and $(P_{n})$ with
corresponding solutions $v_0$ and $u_{n}$),
regularity up to the boundary and  L'H\^opital theorem (see \cite{fmm} for details). 
Therefore, 
$u_n(x)\geq v_0(x)>0$ for all $x\in\Omega$. Thus, by passing to the limit, we conclude that $u$ is a positive
solution of problem $(P)$ and the proof of the theorem is thus complete. \qed


\end{document}